\documentclass[12pt]{amsart}
\usepackage{graphicx,epsfig,amscd,graphics,xypic}
\xyoption{all}

\usepackage[all]{xy}
\usepackage{amssymb}   % For Latex2e
\usepackage{amsmath}
\usepackage{amsthm}

\usepackage{enumerate}
\usepackage{times}
\usepackage[latin1]{inputenc}
\usepackage{psfrag}
\usepackage[active]{srcltx}
\usepackage{color}

\newtheorem{Theorem}{Theorem}[section]
\newtheorem{lemma}[Theorem]{Lemma}
\newtheorem{proposition}[Theorem]{Proposition}
\newtheorem{corollary}[Theorem]{Corollary}

\theoremstyle{definition}
\newtheorem{definition}[Theorem]{Definition}
\newtheorem{remark}[Theorem]{Remark}
\newcommand{\Pic}{\operatorname{Pic}}

\newcommand{\Z}{{\mathbb Z}}

\newcommand{\C}{{\mathbb C}}
\newcommand{\p}{{\mathbb P}}

\newcommand{\codim}{\operatorname{codim}}
\def\leq{\leqslant}
\def\geq{\geqslant}

\def\bibaut#1{{\sc #1}}

\def\phi{\varphi}
\def\ro[#1]{{\textcolor{red}{#1}}}
  
\begin{document}

%\begin{abstract}An embedded manifold is {\it dual defective} if its dual variety is not a hypersurface. Using the geometry of the variety of lines through a general point, we characterize scrolls among dual defective manifolds. This leads to an optimal bound for the dual defect, which improves results due to Ein. We also discuss our conjecture that every dual defective manifold with cyclic Picard group should be secant defective, of a very special type, namely a local quadratic entry locus variety.
%\end{abstract}

%\keywords{Fano manifold, covered by lines, dual and secant defective, scroll}

\title{Remarks on defective Fano manifolds}

\dedicatory{Dedicated to the memory of Alexandru Lascu}
\author{Paltin Ionescu\ and Francesco Russo}
\address{P. Ionescu, Dipartimento di Matematica e Informatica, Universit\` a degli Studi di Ferrara, via Machiavelli, 30, 44121 Ferrara, Italy,   \newline  Email: Paltin.Ionescu\@unife.it }

\address{ F. Russo, Dipartimento di Matematica e Informatica,
Universit\` a degli Studi di Catania,
Viale A. Doria, 6,
95125 Catania, Italy \newline Email: frusso\@dmi.unict.it} 
%\email{Paltin.Ionescu@unife.it}
%\author[Francesco Russo]{Francesco Russo}
%\%address{\sc Dipartimento di Matematica e Informatica\\
%Universit\` a degli Studi di Catania\\
%Viale A. Doria, 6\\
%95125 Catania\\ Italy}
%\email{frusso@dmi.unict.it}

\date{}

\begin{abstract}
 This note continues our previous work on special secant defective (specifically, conic connected and local quadratic entry locus) and dual defective manifolds. These are now well understood, except for the prime Fano ones. Here we add a few remarks on this case, completing the results in our papers \cite{LQEL I}, \cite{LQEL II}, \cite{CC}, \cite{HC} and \cite{DD};
 see also the recent book \cite{Ru}.

\end{abstract}
\maketitle
\section{Introduction}

Let $X\subset {\p}^N_\C:={\p}^N$ be a closed subvariety, that we always assume to be smooth, irreducible, non-degenerate, of dimension $n$ and codimension $c$. The expected dimension of $SX$, the  secant variety of $X$, is $2n+1$. The difference between the expected and the actual dimension is the number $\delta\geq 0$, called the secant defect. If $\delta >0$, $X$ is said to be secant defective; this is the same as being (possibly after an isomorphic projection) of small codimension, in the sense that $n\geq c$. At present we understand better only some subclasses of secant defective manifolds: conic connected (CC for short) and local quadratic entry locus (LQEL for short) ones. CC manifolds are defined by the presence of a conic contained in $X$ and passing through two given general points. An important particular case of CC manifolds are the LQEL ones, defined by the presence of a quadric of dimension $\delta$, contained in $X$ and passing through two given general points; they were introduced in \cite{LQEL I} and further studied in \cite{LQEL II}.

An embedded Fano manifold is prime if its Picard group is generated by the hyperplane section class, $H\in \Pic(X)$. In this case, its index $i$ is defined by the relation $-K_X=iH$ in $\Pic(X)$, $K_X$ being the canonical class. The CC manifolds are Fano and classified, but for the case they are prime \cite{CC}.

In sections 3 and 4, we investigate the property of being CC or LQEL for a Fano manifold of high index.  For instance, we show that for Fano manifolds of index $i$ with $\delta \geq 2$ and $SX=\p^N$,  ``high index'' implies  ``small codimension'', see Corollary~\ref{cor5}~(i). We characterize LQEL manifolds with $\delta \geq 3$ by the equality $i=\frac{n+\delta}{2}$, see Proposition~\ref{prop4}. In Proposition~\ref{prop41} we obtain a classification result for Fano manifolds of high index and small codimension, based on \cite{LQEL I}.

If $X\subset \p^N$ is as before, it is called dual defective (DD for short) if its dual variety $X^*\subset {\p^N}^*$ is not a hypersurface. The number $k=N-1-\dim(X^*)$ is the dual defect of $X$. In section 5 we give ample evidence for our belief that prime Fano DD manifolds are LQEL, improving on \cite{DD}. In particular, we show that a DD manifold is LQEL if (and only if) the relation $\delta=k+2$ holds, Proposition~\ref{prop6}.

In the last section we make some remarks about the case of quadratic manifolds (that includes all known examples of both LQEL and DD prime Fano manifolds in their natural embedding). We also mention the link with the famous Hartshorne Conjecture on manifolds of small codimension, that was shown to hold in the quadratic case \cite{HC}.

We thank the anonymous referee for her/his suggestions that improved both the exposition and point (v) of our 
Proposition~\ref{prop7}.

\section{Preliminaries}

 We work over $\C$ and use the customary notation in algebraic geometry, as in \cite{DD}; in particular, $\p(E)$ is intended in Grothendieck's sense.
We denote by $X\subset \p^N$ a smooth, irreducible, non-degenerate embedded projective variety of dimension $n$ and codimension $c$. For a point $x\in X$, we write ${\bf T}_xX$ for the Zariski tangent space at $x$ and $T_xX$ for its projective closure in $\p^N$. $\mathcal{P}_X$ denotes the first jet bundle of $\mathcal{O}_X(1)$. If $N_{X/\p^N}$ is the normal bundle of $X$ in $\p^N$, we have the standard exact sequence:

\begin{equation}\label{tangent}
0\to N_{X/\p^N}^*(1)\to \mathcal{O}_X^{N+1}\to \mathcal P_X\to 0.
\end{equation}

\begin{definition} $X$ is {\it Fano} if $-K_X$ is ample. $X\subset \p^N$ is called {\it prime Fano} if it is Fano and $\Pic(X)=\Z[H]$, where $H$ denotes the hyperplane section class; it is called {\it Fano of index} $i>0$ if $-K_X=iH$.

\end{definition}

\begin{remark}

Note that our definition of the index of an {\it embedded} Fano manifold is different from the usual one, but is the same if $X$ is prime.

\end{remark}

\begin{lemma}\label{lemma-1}

\begin{enumerate}[(i)]
\item[{\rm (i)}] Let $X$ be a Fano manifold (of dimension $n$) and assume that $-K_X=jD$, for some ample divisor $D$.
Then $j\leq n+1$.
\item[{\rm (ii)}] Let $X\subset \p^N$ be a Fano manifold of index $i\geq \frac{n+3}{2}$. Then $X$ is prime Fano.
\end{enumerate}
 
\end{lemma}

\proof

(i) follows from vanishing results, see \cite{K-O}.

(ii) A result by Wi\'sniewski \cite{Wi} shows that $\Pic(X)$ is cyclic. By (i), $H$ cannot be divisible, so $X$ is prime.

\endproof

Prime Fano manifolds of high index, other than complete intersections, are quite special.

\begin{definition}
$X\subset \p^N$ is {\it covered by lines} if for any point $x\in X$, there is a line $l\subset X$ such that $x\in l$.

\end{definition}

The following fundamental fact, coming from \cite{Mori}, shows that Fano manifolds of high index are covered by lines.

\begin{lemma}\label{lemma0}
If $X\subset \p^N$ is Fano of index $i$ and $i\geq \frac{n+2}{2}$, then $X$ is covered by lines.

\end{lemma}

If $X\subset \p^N$ is covered by lines and $x\in X$ is a general point, we denote by $\mathcal{L}_x\subset \p({{\bf T}}^*_xX)$ the Hilbert scheme of lines contained in $X$ and passing through $x$. $\mathcal{L}_x$ is smooth and  we let $a:=\dim(\mathcal{L}_x)$. If $X\subset \p^N$ is Fano of index $i$, we have $i=a+2$ and $\mathcal{L}_x$ is equidimensional.

\begin{definition}

The {\it secant variety} of $X\subset \p^N$, denoted by $SX$, is the closure of the locus of its secant lines. We have that $\dim(SX)\leq 2n+1$ and $\delta:=2n+1-\dim(SX)$ is the {\it secant defect} of $X$. $X$ is {\it secant defective} if $\delta > 0$. The {\it tangent variety} of $X$, denoted by $TX$, is $TX=\bigcup_{x\in X}T_xX$.
Let us also recall that, if $p\in SX$ is a general point, the {\it entry locus} of $X$ with respect to $p$, denoted by 
$\Sigma_p(X)$, is the intersection of $X$ with the cone of secants to $X$ passing through $p$. The secant defect of $X$ is 
$\dim(\Sigma_p(X))$, for $p\in SX$ a general point.

\end{definition}

\begin{lemma}\label{lemma00}

Let $X\subset \p^N$ be prime Fano of index $i\geq \frac{n+2}{2}$ and let $x\in X$ be a general point. Then one of the following holds:

\begin{enumerate}[(i)]
\item[{\rm (i)}] $S\mathcal{L}_x=\p^{n-1}$, or

\item[{\rm (ii)}]$N_{\mathcal{L}_x/\p^{n-1}}(-1)$ is not ample.

\end{enumerate}
\end{lemma}
\proof
If $X\subset \p^N$ we have that $\delta \geq n-c+1$ and equality holds if and only if $SX=\p^N$. Apply this to $\mathcal{L}_x\subset \p^{n-1}$ (note that $\mathcal{L}_x\subset \p^{n-1}$ is non-degenerate by \cite[Thm. 2.5]{Hwang}). Now, assume that (i) does  not hold. It follows that $\delta(\mathcal{L}_x) > a-(n-1-a)+1\geq 0$. By Terracini's Lemma, see for instance \cite[Thm. 1.4.1]{Ru}, there is a hyperplane in $\p^{n-1}$ such that its tangency locus contains the general entry locus $\Sigma_p(\mathcal{L}_x)$. As $\dim(\Sigma_p(\mathcal{L}_x))=\delta(\mathcal{L}_x) > 0$, (ii) follows.
 
\endproof
\begin{definition}[\cite{LQEL I}]
$X\subset \p^N$ is a {\it local quadratic entry locus} variety, abbreviated LQEL, if, given two general points $x, x' \in X$, there is a quadric of dimension $\delta >0$ contained in $X$ and containing the points $x, x'$.

\end{definition}

Note that any quadric $Q\subset X$ passing through the general points $x, x' \in X$ is part of the entry locus of $X$ with respect to a general point $p$ of the line $\langle x, x' \rangle$. So the dimension of $Q$ is at most the secant defect of $X$.
Thus, LQEL varieties are characterized by the fact that the general entry locus is a (finite) union of quadrics. LQEL varieties are the simplest secant defective manifolds; they were introduced in \cite{LQEL I} and further studied in \cite{LQEL II} and \cite{CC}.

\begin{definition}
$X\subset \p^N$ is {\it conic connected}, abbreviated CC, if, given two general points $x, x'\in X$, there is a conic $C\subset X$ such that $x, x' \in C$.
$X$ is {\it connected by degenerate conics} if $C$ is a pair of incident lines.

\end{definition}

Obviously, LQEL manifolds are CC, and CC manifolds are secant defective.
CC manifolds were studied in \cite{CC}, where a complete classification was obtained, except for the case of prime Fano manifolds. In particular, the same applies for LQEL manifolds.

\begin{definition}
$X\subset \p^N$ is {\it dual defective}, of dual defect $k$, if $\dim(X^*)=N-1-k$, where $X^*\subset {\p^N}^*$ is the dual of $X$ and $k>0$.

\end{definition}

In \cite{BFS} the classification of DD manifolds was reduced to that of prime Fano ones.

Let $X$ be a Fano manifold of dimension $n$ with $-K_X=iH$, $i > 0$, for some ample divisor $H$; let $E$ be a spanned vector bundle on $X$, of rank $r$. 
Let $Y:=\p(E)$, let $M:=\mathcal{O}_Y(1)$ and assume that $M^{n+r-1}=0$. Let $\phi:Y\to W$ be the Stein factorization of the map associated to $|M|$ and let $F$ be its general fiber. Finally, assume that $\det(E)=jH$. 

\begin{proposition}\label{prop1}

Under these assumptions, $F$ is Fano if and only if $j < i$. If this is true and $\Pic(X)$ is cyclic, then $Y$ is also Fano and  $\phi$ is its second Mori contraction.
\end{proposition}
\proof

We have $-K_Y=rM-\pi^*(K_X+\det(E))$, where $\pi:Y\to X$ is the projection. Restricting to $F$ gives: $-K_F=-(K_X+\det(E))|_F=(i-j)H_F$. The rest is clear, since $\Pic(Y)=\Z \oplus \Z$ and the cone of curves of $Y$ is generated by two classes: one of a curve contracted by $\pi$ and one of a curve contracted by $\phi$.

\endproof
\begin{corollary}\label{cor1}

Let $X\subset \p^N$ be Fano of index $i$.

\begin{enumerate}[(i)]
\item[{\rm (i)}] Assume $\delta \geq 2$ and let $E=\mathcal P_X$. Then $F$ is Fano if and only if  $i\geq\frac{n+2}{2}$. In particular $i\geq \frac{n+2}{2}$ implies $i\leq\frac{n+\delta+1}{2}$ with equality holding  if and only if $F\simeq\p^{\delta-1}$.

\item[{\rm (ii)}] If $\delta \geq 2$ and $X$ is prime, then   $\phi:\p(\mathcal P_X)\to W$ is a Mori contraction if and only if $i\geq \frac{n+2}{2}$.

\item[{\rm (iii)}] Assume that $k >0$ and let $E:=N_{X/\p^N}(-1)$. Then $\phi:\p(E)\to W$ is a Mori contraction, $i=\frac{n+k+2}{2}$ and $X$ is prime.
\end{enumerate}
\end{corollary}
\proof

(i) From \eqref{tangent} we have the natural map $\phi': \p(\mathcal {P}_X)\to TX\subset \p^N$, whose Stein factorization is $\phi:\p(\mathcal {P}_X)\to W$. Since $\delta>0$, we have $\dim(TX)=2n+1-\delta$, see \cite[Thm. 1.4]{Zak}. Therefore $\dim(F)=\delta-1$. Reason as in the proof of Proposition~\ref{prop1} and observe that in our case $j=n+1-i$ and $i-j=2i-(n+1)$. If $F$ is Fano, then $2i-(n+1)=i-j\leq \dim(F)+1=\delta$ with equality holding if and only if $F\simeq\p^{\delta-1}$, see \cite{K-O}. Part (ii) is now clear.

(iii) The dual of the exact sequence \eqref{tangent} gives the map defining $X^*$, the dual of $X$, $\phi':\p(N_{X/\p^N}(-1))\to X^*\subset {\p^N}^*$; its Stein factorization is $\phi:\p(N_{X/\p^N}(-1))\to W$, $W$ being the normalization of $X^*$. Reason as above, note that $F$ is a linear space of dimension~$k$ and that $j=n+1-i$.
Since we have $i=\frac{n+k+2}{2}$, we get $i\geq \frac{n+3}{2}$ and  $X$ is prime by Lemma~\ref{lemma-1}~(ii).
%Corollary 2

%Assume $X\subset \p^N$ is Fano, with $Pic(X)=\Z[H]$ ($H$ hyperplane section). If $X$ is DD and $\delta \geq 2$, then $\delta \geq k+2$.

%Proof

%Just put $a=\frac{n+k-2}{2}$ in the inequality from Corollary 1.
\endproof

\begin{remark}
 
\begin{enumerate}[(i)]
\item[{\rm (i)}] For prime Fano manifolds of index $i$, the same condition $i\geq \frac{n+2}{2}$, as in Corollary~\ref{cor1}~(i), ensures that $\mathcal{L}_x \subset \p^{n-1}$ is (non-empty and) non-degenerate, see \cite[Thm. 2.5]{Hwang}.

\item[{\rm (ii)}] Corollary~\ref{cor1}~(iii) proves the Landman Parity Theorem for Fano manifolds: if $k > 0$, $n$ and $k$ have the same parity. A variant of this argument proves the general case, when $X$ is only assumed to be dual defective, see \cite[Prop. 3.1]{DD}.

\end{enumerate}
\end{remark}
\section{Fano manifolds and CC manifolds}

 If $X\subset \p^N$ is CC, then $X$ is Fano and classified unless it is prime of index $i$; if it is prime of index $i$, then $\delta >0$ and $i\geq \frac{n+1}{2}$, see \cite{CC} and \cite[Prop. 3.2]{LQEL II}.

{\it Is it true that, conversely, if $X\subset \p^N$ is a prime Fano manifold of index $i$, such that $i\geq \frac{n+1}{2}$ and $\delta >0$, then $X$ is CC ?}

The question is motivated by the following result of Bonavero--H\" oring, \cite{BH}, which shows that the answer is affirmative in the ``standard" case:

\begin{proposition}[\cite{BH}]\label{prop12}
Let $X\subset \p^N$ be a complete intersection. Assume that $X$ is Fano of index $i\geq \frac{n+1}{2}$. Then $X$ is CC.

\end{proposition}

\begin{lemma}\label{lemma14}
Let $X\subset \p^N$ be a complete intersection. Assume that $X$ is Fano of index $i\geq \frac{n+r}{2}$, for some $r > 0$.
Then $\delta \geq r$.
\end{lemma}
\proof
Let $(d_1,...,d_c)$ be the type of $X$. We have $i=n+1-\sum_{j=1}^c(d_j-1)$, which gives
$2c \leq 2\sum_{j=1}^c(d_j-1) \leq n-r+2$. Finally, we get $\delta \geq n-c+1 \geq r+c-1 \geq r$.

\endproof
\begin{remark}\label{remark3}
\begin{enumerate}[(i)]
\item[{\rm (i)}] The above lemma shows that for a Fano complete intersection of index $i$, the condition $i\geq \frac{n+1}{2}$ implies that $\delta > 0$. 
\item[{\rm (ii)}] The Lagrangian Grassmannian $\mathbb {LG}(2,5)\subset \p^{13}$ is a prime Fano manifold of dimension 6 and index 4 (so we have $i=\frac{n+2}{2}$), but it has $\delta=0$, see the tables in \cite{Ka}.
\end{enumerate}
\end{remark}

At present the answer to the above question is not known and the only (very) partial results require stronger assumptions  and seek for a stronger conclusion, namely that $X$ be connected by degenerate conics. In particular, we need to assume that $\delta \geq 2$.

\begin{lemma}\label{lemma1}

\begin{enumerate}[(i)]

\item[{\rm (i)}] If $X\subset \p^N$ is connected by degenerate conics, then $\delta\geq 2$;

\item[{\rm (ii)}] assume that $X$ is Fano of index $i$ and CC; then $i\geq \frac{n+2}{2}$ if and only if $X$ is connected by degenerate conics. 

\end{enumerate}
\end{lemma}

\proof

(i) Fix two general points $x, x' \in X$. By assumption, there exists a degenerate conic $C_{x,x'}=l\cup l'\subset X$ with $x \in l$ and $x' \in l'$.
 Since $x$ and $x'$ are general points, $N_{l/X}$ and $N_{l'/X}$ are semiample and $C_{x, x'}$ may be deformed inside $X$, by keeping the point $x$ fixed, to a smooth conic, say $C'_{x}$, see \cite[II.7.6.2]{Kol}. Deformations of $C_{x, x'}$ keeping $x$ fixed fill up $X$, so the same is true for the deformations of $C'_{x}$. Therefore, through the general points $x, x''\in X$ there is a smooth conic $C'_{x, x''}$. By hypothesis, there is also a degenerate conic $C''_{x, x''}$ through these points. The union $C'_{x, x''}\cup C''_{x, x''}$ is part of the entry locus $\Sigma_p$ with respect to a general point $p$ of the line $\langle x, x''\rangle$. But $\Sigma_p$ is smooth at the points $x, x''$ and this shows that $\delta=\dim(\Sigma_p)\geq 2$.

(ii) This follows from \cite[Prop. 3.2]{LQEL II}.

\endproof

\begin{proposition}[\cite{H-K}]\label{prop2}
Let $X\subset \p^N$ be Fano of index $i$, with $i\geq \frac{n+3}{2}$. 
Then the following are equivalent:
\begin{enumerate}[(i)]
\item[{\rm (i)}] $S\mathcal{L}_x=\p^{n-1}$ for a general point $x\in X$;  

\item[{\rm (ii)}] $X$ is connected by degenerate conics.

\end{enumerate}

\end{proposition}

 In the sequel we only need that (i) implies (ii), so we include a very simple proof of this result, based on the Terracini Lemma.

\proof

(i) implies (ii). We remark first that Lemma~\ref{lemma-1} implies that $X$ is prime Fano; moreover, by Lemma~\ref{lemma0} $X$ is covered by lines, so (i) makes sense. The condition $i\geq \frac{n+3}{2}$ is equivalent to $a\geq \frac{n-1}{2}$; therefore $\mathcal{L}_x\subset \p^{n-1}$ is smooth, irreducible and non-degenerate, for $x\in X$ a general point, see \cite{Hwang}.
Let $l\subset X$ be a general line and let $x, x' \in l$ be general points. 
Let $[l']\in \mathcal{L}_{x'}$ be a general line and let $e\in l'$ be a general point. Denote by $C(x)$ the locus of lines through $x\in X$; it is a cone of dimension $a+1$.

By the Terracini Lemma, the condition $S\mathcal{L}_x=\p^{n-1}$ is equivalent to 
\begin{equation}\dim(T_{x}C(x')\cap T_{e}C(x'))=2a+2-n.\label{unu}\end{equation}
Denote by E the locus of points lying on lines meeting $C(x)$. $X$ is connected by degenerate conics exactly when $E=X$; equivalently, $\dim(E) \geq n$.

If $\pi:\mathcal{Y}\to \mathcal{L}$ is the family of lines covering $X$ and $\psi:\mathcal{Y}\to X$ is the projection, let $V:=\pi^{-1}(\pi(\psi^{-1}(C(x))))$.    Note that $\dim(V)=2a+2$. $\dim(E) \geq n$ is the same as asking that, for some $e\in V$, the fiber $V_e$ of the map
$\psi|_V:V\to E$ satisfies $\dim(V_e)\leq 2a+2-n$. But $\dim(V_e)=\dim_{x'} (C(x)\cap C(e))$.

Therefore, keeping the above notation, to prove that $X$ is connected by degenerate conics, we have to show that 

\begin{equation}\dim_{x'}(C(x)\cap C(e))\leq 2a+2-n.\label{!}\end{equation}
Let us recall from \cite[Prop. 2.2]{DD} that, if $x, x'$ are general points of a general line $l\subset X$, we have 
\begin{equation}T_xC(x')=\bigcap_{y\in l}T_yX=T_{x'}C(x).\label{doi}\end{equation}
Now, using \eqref{unu} and \eqref{doi}, we can write
\begin{align*}\dim_{x'}(C(x)\cap C(e))&\leq \dim(T_{x'}C(x)\cap T_{x'}C(e))\\&=\dim(T_xC(x')\cap T_eC(x'))=2a+2-n\end{align*} and \eqref{!} is proved.\endproof

Note that, since $X$ is smooth, the inequality $\dim_{x'}(C(x)\cap C(e))\geq 2a+2-n$ also holds; therefore, the cones $C(x)$ and $C(e)$ intersect transversally
at $x'$.
\endproof

\begin{corollary}[\cite{H-K}, {\rm see also} \cite{CC}] \label{cor2}
Let $X\subset \p^N$ be Fano of index $i$ and assume that $i > \frac{2n}{3}$. Then $X$ is connected by degenerate conics.

\end{corollary}

Remark \ref{remark3}~(ii) shows that the bound $i > \frac{2n}{3}$ is optimal.

Proposition \ref{prop2} implies Proposition \ref{prop12} in case $i\geq \frac{n+2}{2}$:

\begin{proposition}\label{prop23}

Let $X\subset \p^N$ be a complete intersection and assume $X$ is Fano, of index $i$. The following are equivalent:

\begin{enumerate}[(i)]
\item[{\rm (i)}] $i\geq \frac{n+2}{2}$;

\item[{\rm (ii)}] $X$ is connected by degenerate conics.

\end{enumerate}
\end{proposition}
\proof
(i) implies (ii):

The result is clear for $n\leq 2$, so we may assume $n\geq 3$ from now on. It follows from (i) that, for $x\in X$ a general point, $\mathcal{L}_x \subset \p^{n-1}$ is a smooth, non-degenerate, complete intersection of positive dimension $a=i-2$. In particular, it is connected and hence irreducible. Moreover, since it is a complete intersection, it follows that $S\mathcal{L}_x=\p^{n-1}$, see Lemma~\ref{lemma00}. Then (ii) follows by the proof of Proposition~\ref{prop2}.

(ii) implies (i) follows from Lemma~\ref{lemma1}.

\endproof
\begin{proposition}\label{prop3}

Let $X\subset \p^N$ be Fano of index $i$ and suppose $\delta \geq 2$. Assume
that 
$i > \frac{n+1}{2}+\frac{\delta}{4}$. Then $X$ is connected by degenerate conics.
\end{proposition}
\proof

Note that Lemma~\ref{lemma-1} implies that $X$ is prime Fano. Let $Y=\p(\mathcal{P}_X)$ and let $\phi:Y\to W$ be as in the proof of Corollary~\ref{cor1}. By Corollary~\ref{cor1}~(i), the general fiber $F$ of $\phi$ is Fano and $i(F)=2i-n-1 > \frac{\dim(F)+1}{2}=\frac{\delta}{2}$. This implies that $F$ is covered by lines, see Lemma~\ref{lemma0}. But $F$ passes through the general point of $Y$; so $Y$  is covered by a family of lines, projecting onto lines of $X$ and contracted by 
the map $\phi':Y\to TX$. We claim that this implies $S\mathcal{L}_x=\p^{n-1}$. By Proposition~\ref{prop2}, it follows that $X$ is conic-connected (by degenerate conics). 

Fix a point $y\in Y$ and let $x=\pi(y)$. Assume that a line $l\subset Y, y\in l$ is contracted to the general point $u\in TX$. It follows that all (projective) tangent spaces to $X$ at points of $\pi(l)$ pass through the point $u$. But, for $s\in l$ a general point, the intersection $\bigcap_{t\in l}T_tX$ is the projective tangent space to the cone of lines through $s$, at a general point of $l$, see \eqref{doi}. As the point $u\in T_xX$ is general, this means that $T\mathcal{L}_x=\p^{n-1}$, so we have a fortiori that $S\mathcal{L}_x=\p^{n-1}$.

\endproof

\begin{corollary}\label{cor5}

Let $X\subset \p^N$ be Fano of index $i$  with $\delta \geq 2$. Then  we have: 
\begin{enumerate}[(i)]
\item[{\rm (i)}] $i\leq \frac{n+\delta}{2}$; in particular, if $i\geq \frac{n+r}{2}$ for some $r > 0$, then $\delta \geq r$ (cf. Lemma \ref{lemma14}). If $SX=\p^N$, this becomes $c\leq n-r+1$.

\item[{\rm (ii)}] The general fiber of the map $\phi:\p(\mathcal{P}_X)\to W$ is not linear.
\end{enumerate}
\end{corollary}

\proof (i) Suppose that $i> \frac{n+\delta}{2}$; since $\delta\geq 2$, we get $i> \frac{n+1}{2}+\frac{\delta}{4}$. By Proposition~\ref{prop3}, $X$ is conic-connected. It follows that $i\leq \frac{n+\delta}{2}$, by \cite[Prop. 3.2]{LQEL II}.

(ii) It follows from the previous inequality and part (i) of Corollary 1.
%Assume that the general fiber, $F$, was linear. We get that $i(F)=\dim(F)+1=\delta$. But the formula from Proposition~\ref{prop1} gives $i(F)=2i(X)-%n-1$. It follows that $2i(X)=n+1+\delta$, which contradicts (i).

\endproof

%If $\delta=2$, the fiber $F$ is a conic; it follows that the irreducible components of the entry locus are quadrics, so $X$ is LQEL. This follows from the remark that the polar with respect to a point $p$ commutes with taking general hyperplane sections through 
%$p$, $p(X)\cap H=p(X\cap H)$.

%Corollary 4

%Assume $X\subset \p^N$ to be Fano, with $Pic(X)=\Z[H]$ and $\delta \geq 2$. If $i\geq \frac{n+s}{2}$, then $s\leq \delta$.

%Question: The above result is clearly false for $\delta=0$; what happens for $\delta=1$?

\section{Fano manifolds and LQEL manifolds}

\begin{proposition}\label{prop4}

Assume that $X\subset \p^N$ is Fano of index $i$ and $\delta \geq 3$. The following conditions are equivalent:
\begin{enumerate}[(i)]
\item[{\rm (i)}] $X$ is LQEL;

\item[{\rm (ii)}] $i=\frac{n+\delta}{2}$;

\item[{\rm (iii)}] the general fiber of the map $\phi:\p(\mathcal{P}_X)\to W$ is isomorphic to the quadric $Q^{\delta-1}$; 

\item[{\rm (iv)}] if $x\in X$ is a general point, we have $S\mathcal{L}_x=\p^{n-1}$ and $\delta(X)\leq \delta(\mathcal{L}_x)+2$.

\end{enumerate}

\end{proposition}
\proof

Since $\delta \geq 3$, $X$ is prime Fano, see \cite{B-L}. (i) implies (ii) and (i) implies (iv) follow from \cite{LQEL I}. If $F$ is a general fiber of $\phi$, the equivalence of (ii) and (iii) comes from the formula in Proposition~\ref{prop1}: $i(F)=2i(X)-n-1$, and the fact that $Q^r$ is characterized by having its index equal to $r$, see \cite{K-O}. (ii) implies (i) comes from Proposition~\ref{prop3} and \cite[Prop. 3.2]{LQEL II}. Finally, from (iv) we infer 
$\delta(X)\leq \delta(\mathcal{L}_x)+2=a-(n-1-a)+1+2=2a+4-n=2i-n$, so (ii) holds by using also Corollary~\ref{cor5}.

\endproof

The next result shows that Fano manifolds of high index and small codimension are very special.

\begin{proposition}\label{prop41}
Let $X\subset \p^N$ be Fano of index $i$. Assume that $i\geq n-\frac{c-1}{2}$ and $n\geq 2c-2$. Then $X$ is an LQEL manifold and one of the following holds (up to an automorphism of $\p^N$):

\begin{enumerate}[(i)]
\item[{\rm (i)}] $c=1, n\geq 1$ and $X$ is $Q^n\subset \p^{n+1}$;

\item[{\rm (ii)}] $n=2c$ and $X$ is one of: $\mathbb G(1, 4)\subset \p^9$ or $S^{10}\subset \p^{15}$, in their natural embeddings;

\item[{\rm (iii)}] $n=2c-1$ and $X$ is one of: $\mathbb G(1, 4)\cap H \subset \p^8$ or $S^{10}\cap H \subset \p^{14}$, where $H, H'$ are general hyperplanes in $\p^N$;

\item[{\rm (iv)}]  $n=2c-2$ and $X$ is one of:  $\mathbb G(1, 4)\cap H\cap H' \subset \p^7$, $S^{10}\cap H\cap H' \subset \p^{13}$, or one of the following Severi varieties: 

$\p^2\times \p^2\subset \p^7$, projection of the Segre embedding into $\p^8$, 
 
$\mathbb G(1,5)\subset \p^{13}$, projection of the Pl\" ucker embedding into $\p^{14}$ and

$E_6\subset \p^{25}$, projection of the natural embedding into $\p^{26}$.

\end{enumerate}
\end{proposition}
\proof
Since $n\geq 2c-2$, it follows from Zak's Linear Normality Theorem, see \cite{Zak}, that we have $SX=\p^N$. 
Assume first that $c\geq 4$; we get $\delta=n-c+1\geq c-1\geq 3$. The hypothesis  $i\geq n-\frac{c-1}{2}$ reads $i\geq \frac{n+\delta}{2}$ and from Corollary \ref{cor5} and Proposition \ref{prop4} it follows that $X$ is an LQEL manifold. Now the classification follows from the results in section 3 of \cite{LQEL I}; due to \cite{Mok}, the last case of \cite[Cor. 3.2]{LQEL I}  is ruled out.

Assume next that $c\leq 3$. For $c=1$ we get case (i) by \cite{K-O}. The case $c=2$ is impossible. Finally, suppose that $c=3$. We infer that $X$ is a del Pezzo manifold, see \cite{Fuj} for their classification. From the condition $c=3$ we deduce $n\geq 4$ and we get the Grassmannian $\mathbb G(1,4)\subset \p^9$ and its linear sections that are mentioned in the cases 
(ii), (iii) and (iv).

\endproof

\begin{proposition}\label{prop5}

Let $X\subset \p^N$ be an LQEL manifold and assume that $\delta \geq 3$. Let $Q$ be a general fiber of  $\phi:\p(\mathcal{P}_X)\to W$. Then we have:
\begin{enumerate}[(i)]
\item[{\rm (i)}] $N_{Q/X}$ is spanned;

\item[{\rm (ii)}] $\mathcal{P}_X|_Q$ is spanned, uniform, with splitting type $(0,\ldots,0,1,\ldots,1)$ and degree $n-1-a$ on any line contained in the quadric $Q$.
\end{enumerate}
\end{proposition}

\proof

(i) The standard diagram 

\begin{equation*}\xymatrix{& &0\ar[d] &0\ar[d]& \\
& & {T_{Y/X}}_{|Q} \ar@{=}[r] \ar[d]  &{T_{Y/X}}_{|Q} \ar[d]&  \\
0\ar[r] &T_Q \ar@{=}[d] \ar[r] & {T_Y}_{|Q}\ar[d] \ar[r] & N_{Q/Y}\ar[r] \ar[d] & 0\\
0\ar[r] &T_Q  \ar[r] & {T_X}_{|Q}\ar[d] \ar[r] & N_{Q/X}\ar[r] \ar[d] & 0\\
& &0 &0&}\end{equation*}
shows that $N_{Q/X}$ is a quotient of the trivial vector bundle $N_{Q/Y}$, where $Y=\p(\mathcal{P}_X)$, see Corollary~\ref{cor1}~(i).

(ii) If $l\subset Q$ is a line, from the standard exact sequence
\begin{equation*}\xymatrix{0\ar[r] & N_{l/Q} \ar[r] & N_{l/X} \ar[r] & N_{Q/X}|_{l}\ar[r] & 0
}\end{equation*}
it follows that $N_{l/X}$ is spanned. Therefore its splitting type is $(0,\ldots,0,1,\ldots,1)$ and the conclusion about $\mathcal{P}_X|_l$ follows from \cite[Prop. 2.2]{DD}.

\endproof

\section{Fano manifolds and DD manifolds}

In this section we collect a number of facts that support our belief that any DD manifold which is Fano of index $i$ is LQEL.

\begin{remark}
\begin{enumerate}[(i)]

\item[{\rm (i)}] If $X\subset \p^N$ is DD, it follows from Zak's Theorem on Tangencies, \cite{Zak}, that we have 
$k \leq c-1$.

\item[{\rm (ii)}] If $X\subset \p^N$ is LQEL and $X\neq Q^n$, we have $\delta \leq c+1$, see \cite[Prop. 4.2]{DD}.
\end{enumerate}
\end{remark}

\begin{proposition}\label{prop6}

Let $X\subset \p^N$ be Fano. Assume that $X$ is DD. Then $X$ is LQEL if and only if $\delta=k+2$. \end{proposition}

\proof

Since $X$ is DD, we have $i=\frac{n+k+2}{2}$; since $X$ is LQEL, we have $i=\frac{n+\delta}{2}$. So, $\delta=k+2$. Conversely, if $k=\delta - 2$, it follows that
$i=\frac{n+\delta}{2}$ and $X$ is LQEL by Proposition~\ref{prop4}.

\endproof

\begin{proposition}\label{prop7}

Let $X\subset \p^N$ be a DD manifold which is Fano of index $i$. Then we have:
\begin{enumerate}[(i)]
\item[{\rm (i)}] If $\delta\geq 2$, then $\delta\geq k+2$; if $\delta=3$, then $X$ is LQEL;

\item[{\rm (ii)}] assume that $k\geq 2$. Then one of the following holds:

\begin{enumerate}[(a)]

\item[{\rm (a)}] $\delta\geq k+2$, or

\item[{\rm (b)}] $n\leq 2c-1$ and for $x\in X$ general, $\mathcal{L}_x\subset \p^{n-1}$ is prime Fano and covered by lines;
\end{enumerate}

\item[{\rm (iii)}] cf. \cite{LQEL II}  if $k > \frac{n-6}{3}$, then $X$ is CC and $\delta\geq k+2$;

\item[{\rm (iv)}] \cite{DD} we have that $k\leq \frac{n+2}{3}$ and equality holds exactly when $X$ is $S^{10}\subset \p^{15}$;

\item[{\rm (v)}] if $k\geq \min\{\frac{n-4}{2}, \frac{n+2}{3}\}$, then $X$ is LQEL;

\item[{\rm (vi)}] cf. \cite{Ein} if $k=c-1$ and $n\leq \max\{2c+2, 3c-5\}$, then $X$ is one of $\mathbb G(1,4)\subset \p^9$ or $S^{10}\subset \p^{15}$; 
%or $\p^1\times \p^{n-1}$;

\item[{\rm (vii)}] cf. \cite{Mu} if $k=c-2$ and $n\leq 2c$, then $X\subset \p^N$ is a hyperplane section of one of the two varieties in (vi).

Assume now that X is LQEL and $X\neq Q^n$; then we have:

\item[{\rm (viii)}] If $\delta=c+1$, then $X\subset \p^N$ is one of the two varieties in (vi);

\item[{\rm (ix)}] cf. \cite{Fu}  $\delta \leq \frac{n+8}{3}$ and equality holds exactly when $X$ is one of $S^{10}\subset \p^{15}$, $E_6\subset \p^{26}$, or its isomorphic projection to $\p^{25}$.
\end{enumerate}\end{proposition}
\proof
By Corollary~\ref{cor1}~(iii), $X$ is prime.

(i) The first part follows from Corollary~\ref{cor5}~(i), given that $i=\frac{n+k+2}{2}$. If $\delta=3$, we get $3=\delta\geq k+2$, so that $k=1=\delta -2$ and Proposition~\ref{prop6} applies.

(ii) Assume first that $n\geq 2c$. It follows that $\delta\geq n-c+1\geq 2$. Therefore $\delta\geq k+2$, by (i). Assume from now on that $n\leq 2c-1$. As we have $i\geq \frac{n+3}{2}$, we know that $\mathcal{L}_x\subset \p^{n-1}$ is smooth, non-degenerate and irreducible, see \cite{Hwang}. If $S\mathcal{L}_x=\p^{n-1}$, by 
Proposition~\ref{prop2} $X$ is connected by degenerate conics, so $\delta\geq 2$, see Lemma~\ref{lemma1}~(i). So we again get that $\delta\geq k+2$ by (i). Now
we may assume that $S\mathcal{L}_x\neq \p^{n-1}$ and hence $\delta(\mathcal{L}_x) > a-(n-1-a)+1=k\geq 2$. By the Barth-Larsen Theorem, \cite{B-L}, we get that $\Pic(\mathcal{L}_x)$ is cyclic, generated by the hyperplane section class. On the other hand, we have $n-2c+2\leq 1\leq k$ and it follows from \cite[Prop. 3.2]{DD} that all lines in $X$ are contact lines. Therefore $\mathcal{L}_x$ is covered by lines. As its Picard group is cyclic, generated by the hyperplane section, it is  prime Fano, of some index $j$.  

(iii) Since we have $2i=n+k+2$, the condition $k > \frac{n-6}{3}$ is equivalent to $i > \frac{2n}{3}$. So $X$ is connected by degenerate conics by  Corollary~\ref{cor2}.  From Lemma~\ref{lemma1}, we have $\delta\geq 2$, so (i) applies.

(iv) See \cite[Cor. 3.5]{DD}.

(v) Using (iv) and \cite{Ein2}, we are left with the cases $n\leq 15$ and $n=2k+4$ or $n=2k+3$. When $n\leq 10$ we may use the classification in \cite{BFS}, together with \cite{Fuj} and \cite{Muk}. The only remaining case is $n=12$, $k=4$, $i=9$.  
By Zak's Theorem on Tangencies we have $k\leq c-1$, so that $c\geq 5$. Therefore, $4=k\geq n-2c+2$, so that all lines on $X$ are contact, see \cite[Prop. 3.2]{DD}. Let $x\in X$ be a general point and $X':=\mathcal{L}_x \subset \p^{11}$. We have $\dim(X')=i-2=7$ and $\codim(X')=4$,
so $X'\subset \p^{11}$ is prime Fano and covered by linear spaces of dimension three. Let $x'\in X'$ be a general point and let $X'':=\mathcal{L}_{x'} \subset \p^{6}$. If $\dim(X'')\geq 3$, we have  $i(X')=\dim(X'')+2\geq 3+2=5$. So, $X'$ would be a del Pezzo or Mukai prime Fano manifold, which is excluded by their classification, see \cite{Fuj} and \cite{Muk}.
Since $X'$ is covered by linear spaces of dimension three, we are left with the case where  $\dim(X'')=2$. The following argument has been suggested by the referee. We let $L\subset X'$ be a linear space (of dimension three) passing through the general point $x'\in X'$ and we let $x'\in l\subset L$ be a line. From the standard exact sequence
\begin{equation*}\xymatrix{0\ar[r] & N_{l/L} \ar[r] & N_{l/X'} \ar[r] & N_{L/X'}|_{l}\ar[r] & 0
}\end{equation*}
it follows that $N_{L/X'}|_{l}$ is nef and its degree is zero, so it is trivial. Applying \cite[Thm. 3.2.1]{OSS}, we get that $N_{L/X'}$ is trivial. Now we may quote \cite[Prop. 6.4]{NO} to exclude this last possibility.
%We may assume that $SX=\p^N$, so $\delta=n-c+1$. It follows that $c\geq k+1=5$, hence $\delta\leq 8$. If $c=5$, we have $\delta=8$.  Otherwise, any general line $l$ is contained in infinitely many $4$-dimensional planes contained in $X$. This is seen by observing that $\dim(\langle \bigcup_{y\in l} T_yX\rangle)\leq \min\{N, 2n-1-a\}=16 < N-1$, cf. \cite{DD}. 

%Therefore, $X'=:\mathcal{L}_x\subset \p^{11}$ is prime Fano, of dimension $7$ and covered by lines. If $x'\in X'$ is a general point, the family of lines passing through $x'$ and contained in $X'$ is of dimension at least $3$. It follows that $X'$ must be a del Pezzo, see \cite{Fuj}, or a Mukai, see \cite{Muk}, manifold and this is impossible by their classification. The case $c=5$, $n=12=2c+2$ remains open. Of course, if it exists, it will contradict the Hartshorne Conjecture !

(vi) If $n\leq 3c-5$, we get $k\geq \frac{n+2}{3}$ and (iv) applies. If $n\leq 2c+2$, we have $k\geq \frac{n-4}{2}$ and it follows from (v) that $X$ is LQEL.  But an LQEL manifold satisfies the condition $\delta=k+2=c+1$. From Lemma~\ref{lemma2} below, the condition $\delta=c+1$ is equivalent to $n=2c$.
The result follows now from \cite[Cor. 3.1]{LQEL I}.

(vii) As above, it follows from (iv) and (v) that $X$ is LQEL. Therefore we have $\delta=k+2=c$, so $c=\delta\geq n-c+1$ and $n\leq 2c-1$. This gives $\delta=c > \frac{n}{2}$ and the result again follows from \cite[Cor. 3.1]{LQEL I}.
%From Lemma~\ref{lemma2} below, we have that $\delta \leq c+1$. If $n\leq c=k+2$, we get that $k\geq n-2$, so $X$ is a scroll. But then $X$ is a linear space, since it is prime Fano. Thus we may assume that $n\geq c+1$, so that $\delta\geq 2$ and, by (i), $\delta \geq k+2=c$. Therefore we are left with two cases:
%$\delta=c=k+2$ or $\delta=c+1=k+3$. In the first case, $X$ is LQEL by Proposition~\ref{prop6} and $\delta=c\geq n-c+1$ gives $\delta > \frac{n}{2}$. The result follows from \cite{LQEL I}. Assume now that $\delta=c+1=k+3$ and $n=2c$. Since $n$ and $k$ have the same parity, $k$ is even, so $c$ is even too. So, we have $k\geq 2$ and $c\geq 4$. The cases $c=4$, $n=8$, $k=2$, or $c=6$, $n=12$, $k=4$, or $c=8$, $n=16$, $k=6$ are excluded, see (v), so we may assume $c\geq 10$. Let $x\in X$ be a general point and $X':=\mathcal{L}_x \subset \p^{2c-1}$. We have $\dim(X')=\frac{n+k-2}{2}=\frac{3c-4}{2}$. Now, we find that $\dim (X')\geq \codim(X') +2$ and $2=k\geq n-2c+2=2$; therefore $X'\subset \p^{2c-1}$ is prime Fano and covered by lines, see \cite{DD}. In fact, $X'$ is covered by linear spaces of dimension $k-1=c-3 > \frac{3c-4}{4}$, for 
%$c\geq 10$. By the results in \cite{Sato} (see also \cite{BI} for a simple proof in the spirit of the present paper), $X'\subset \p^{2c-1}$ must be a scroll. But then $X'$ must be a linear space, since it is prime. This would imply that $X$ is linear. 

(viii) From Lemma~\ref{lemma2} below, $\delta=c+1$ is equivalent to $n=2c$, so the result follows from \cite[Cor. 3.1]{LQEL I}.

(ix) We may assume that $\delta\geq 5$. We know that $X':=\mathcal{L}_x\subset \p^{n-1}$ is also an LQEL, having dimension $a=\frac{n+\delta-4}{2}$ and secant defect $\delta'=\delta-2$, see \cite[Thm. 2.8]{LQEL I}. The bound $\delta \leq \frac{n+8}{3}$ follows quickly from this, see \cite[Cor. 4.4]{DD}. If equality holds, we get that $\delta' > \frac{a}{2}$ and from \cite[Cor. 3.1]{LQEL I} we deduce that $X'$ is one of $\mathbb G(1,4)\subset \p^9$ or $S^{10}\subset \p^{15}$. Using \cite{Mok}, we get the three cases announced in (ix).

\endproof

\begin{lemma}\label{lemma2}

Let $X\subset \p^N$ be smooth, non-degenerate. The following conditions are equivalent:
\begin{enumerate}[(i)]
\item[{\rm (i)}] $n\leq 2c$;

\item[{\rm (ii)}] $\delta\leq c+1$.
\end{enumerate}
Moreover, equality holds in {\rm (i)} if and only if it holds in {\rm (ii)}.
\end{lemma}
\proof

(ii) implies (i):
We have $c+1\geq \delta\geq n-c+1$, so $n\leq 2c$.

(i) implies (ii): 
If $SX=\p^N$, we have $\delta=n-c+1$ and (ii) is clear. Otherwise, by Zak's Linear Normality Theorem, see \cite{Zak}, we have $\delta\leq \frac{n}{2}$, so $\delta\leq c$.

The other assertion of the lemma is proved in a similar way.
\endproof

\begin{proposition}\label{prop8}

Assume that $X$ is Fano and DD. Let $\p^k\cong L\subset X$ be the contact locus of a general hyperplane. Then we have:
\begin{enumerate}[(i)]

\item[{\rm (i)}] $N_{L/X}$ is spanned;

\item[{\rm (ii)}] $\mathcal{P}_X|_L$ is spanned, uniform, with splitting type $(0,\ldots,0,1,\ldots,1)$ and degree $n-1-a$ on any line contained in $L$.
\end{enumerate}\end{proposition}
\proof

We proceed as in the proof of Proposition~\ref{prop5}. Let $Y:=\p(N_{X/\p^N}(-1))$ and note that $L$ is the general fiber of the morphism $\phi:\p(N_{X/\p^N}(-1))\to W$, see Corollary~\ref{cor1}~(iii). Therefore the normal bundle $N_{L/Y}$ is trivial and a diagram similar to the one from the proof of Proposition~\ref{prop5} shows that $N_{L/X}$ is spanned. Now, if $l\subset L$ is a line, the exact sequence:
\begin{equation*}\xymatrix{0\ar[r] & N_{l/L} \ar[r] & N_{l/X} \ar[r] & N_{L/X}|_{l}\ar[r] & 0
}\end{equation*}
gives that $N_{l/X}$ is spanned and the conclusions follow from \cite[Prop. 2.2]{DD}.

\endproof

\section{Remarks on the quadratic case and the Hartshorne Conjecture}

\begin{definition}
$X\subset \p^N$ is {\it quadratic} if it is scheme-theoretically an intersection of quadrics.
\end{definition}

\begin{remark}
\begin{enumerate}[(i)]
\item[{\rm (i)}] All known examples of linearly normal LQEL or DD prime Fano manifolds are quadratic.

\item[{\rm (ii)}] All known examples of prime Fano manifolds of high index, in their natural linearly normal embedding,  are either complete intersections or quadratic.

\item[{\rm (iii)}] Quadratic manifolds of small codimension ($n\geq c+1$) are quite restricted: they are Fano, covered by lines and $\mathcal{L}_x \subset \p^{n-1}$ is also quadratic,  see \cite{HC}.

\item[{\rm (iv)}] If $X\subset \p^N$ is quadratic and covered by lines, the splitting-type of $N_{X/{\p^N}}(-1)|l$ is 
$(0,...,0,1,...,1)$ if $l$ is a general line. We recall from \cite[Prop. 3.2]{DD} that, if $X$ is DD and $k\geq n-2c+2$, this property still holds.

\end{enumerate}
\end{remark}

Let us also recall the famous {\it Hartshorne Conjecture}, HC for short: 
\smallskip

{\it If $n\geq 2c+1$, $X$ is a complete intersection.} 

\smallskip
\noindent We consider it very plausible for (prime) Fano manifolds. The HC holds in the following special cases: for Fano manifolds in codimension two \cite{B-C}, for quadratic manifolds \cite{HC}, for LQEL manifolds and for $\mathcal{L}_x \subset \p^{n-1}$ if $X$ is DD, see \cite{DD}. 

Both LQEL and DD Fano manifolds have high index and we expect the HC to be easier to prove in this case. An example is Proposition~\ref{prop41}, where ``high index'' means $i\geq n-\frac{c-1}{2}$. This inequality excludes all complete intersections, but for the hyperquadrics. We close this note with another partial result of the same flavor, where only complete intersections should be present.

\begin{proposition}
Let  $X \subset \p^N$ be Fano of index $i\geq \frac{N+2}{2}$. Then the following hold:
\begin{enumerate}[(i)]
\item[{\rm (i)}] $X$ is prime and covered by lines; 

\item[{\rm (ii)}] $X$ is connected by degenerate  conics;

\item[{\rm (iii)}] $SX=\p^N$;

\item[{\rm (iv)}] $n\geq 2c+1$.
\end{enumerate}
\end{proposition}

\proof
(i) The fact that $X$ is prime follows from Lemma~\ref{lemma-1}; Lemma~\ref{lemma0} shows that $X$ is covered by lines.

(ii) If $x, x' \in X$ are general points, let $C(x), C(x')$ be the cones of lines passing through the respective points and contained in $X$; they have dimension $a+1$. Since we have $i=a+2$, it follows that $2a+2\geq N$, so the cones $C(x)$ and $C(x')$ must intersect.

(iii) and (iv) From (ii) and Lemma~\ref{lemma1} it follows that $\delta \geq 2$. Now, Corollary~\ref{cor5} yields that 
$\delta \geq c+2$. By Lemma~\ref{lemma2} we get that $n\geq 2c+1$ and the conclusion follows from Zak's Linear Normality Theorem.

\endproof

Properties (i), (ii) and (iii) are necessary conditions for $X$ to be a complete intersection and, by (iv), the HC predicts this is indeed the case. Note that the inequality $i\geq \frac{N+2}{2}$ is optimal, as the examples $\mathbb G(1,4)\subset \p^9$ and $S^{10}\subset \p^{15}$ show.


\begin{thebibliography}{Comemee}


%\bibitem[1]{A}  \bibaut{C. Araujo}, Rational curves of minimal degree and characterizations of projective spaces, {\it Math. Ann.} {\bf 335} (2006), 937--951.
 
\bibitem[1]{B-C} \bibaut{E. Ballico, L. Chiantini}, On smooth subcanonical varieties of codimension $2$ in $\p^n$, 
$n\geq 4$, {\it Ann. Mat. Pura Appl.} {\bf 135} (1983), 99--117.

\bibitem[2]{B-L} \bibaut{ W. Barth, M.E. Larsen},  On the homotopy
groups of complex projective algebraic manifolds, {\it Math. Scand.} {\bf  30} (1972), 88--94.

\bibitem[3]{BFS}  \bibaut{M.C. Beltrametti, M.L. Fania, A.J. Sommese}, On the discriminant variety of a projective manifold, {\it Forum Math.} {\bf 4} (1992), 529--547.

\bibitem[4]{BH} \bibaut{L. Bonavero, A. H\" oring}, Counting conics in complete intersections, 
{\it Acta Math. Vietnamica} {\bf 35} (2010), 23--30.

%\bibitem[BI] {BI} Beltrametti-Ionescu

%\bibitem[Deb]{Deb} Debarre



\bibitem[5]{Ein} \bibaut{L. Ein}, Varieties with small dual variety. I, {\it Invent.  Math.}   {\bf 86} (1986), 63--74.


\bibitem[6]{Ein2} \bibaut{L. Ein}, Varieties with small dual variety. II, {\it Duke  Math.  J.}  {\bf 52} (1985), 895--907.

\bibitem[7]{Fu} \bibaut{B. Fu}, Inductive characterizations of hyperquadrics, {\it Math. Ann.} {\bf 340} (2008), 185--194.


\bibitem[8]{Fuj} \bibaut{T. Fujita}, {\it Classification Theories of Polarized Varieties}, London Math. Soc. Lecture Note Ser., vol. 155, Cambridge Univ. Press, 1990.

\bibitem[9]{Hwang}  \bibaut{J.M. Hwang}, Geometry of minimal rational curves on Fano manifolds, in {\it School on Vanishing Theorems and Effective Results in Algebraic Geometry (Trieste, 2000)}, ICTP Lect.\ Notes, vol.\ 6, Abdus Salam Int.\ Cent.\ Theoret.\ Phys., 2001, pp.\ 335--393.

\bibitem[10]{H-K} \bibaut{J.M. Hwang, S. Kebekus}, Geometry of chains of minimal rational curves, {\it J. Reine Angew. Math.} {\bf 584} (2005), 173--194.

\bibitem[11]{LQEL II} \bibaut{P. Ionescu, F. Russo},  Varieties with quadratic entry locus. II, {\it Compositio Math.}  {\bf 144} (2008), 949--962.

\bibitem[12]{CC} \bibaut{P. Ionescu, F. Russo}, Conic-connected manifolds,  {\it J. Reine Angew. Math.} {\bf 644} (2010), 145--157.

\bibitem[13]{HC}  \bibaut{P. Ionescu, F. Russo}, Manifolds covered by lines and the Hartshorne Conjecture for quadratic manifolds, {\it Amer. J. Math.} {\bf 135} (2013), 349--360.


\bibitem[14]{DD} \bibaut{P. Ionescu, F. Russo}, On dual defective manifolds, {\it Math. Res. Lett.} {\bf 21} (2014), 
1137--1154.

\bibitem[15]{Ka} \bibaut{H. Kaji}, Homogeneous projective varieties with degenerate secants, {\it Trans. Amer. Math. Soc.} {\bf 351} (1999), 533--545.

\bibitem [16]{K-O} \bibaut{S. Kobayashi, T. Ochiai}, Characterizations of complex projective spaces and hyperquadrics, {\it J. Math. Kyoto Univ.} {\bf 13} (1973), 31--47. 

\bibitem [17]{Kol} \bibaut{J. Koll\'ar}, {\it Rational Curves on Algebraic Varieties}, Ergeb. Math. Grenzgeb. (3), vol. 32, Springer, 1996.

\bibitem[18]{Mok} \bibaut{N. Mok}, Recognizing certain rational homogeneous manifolds of Picard number 1 from their varieties of minimal rational tangents, in {\it Third International Congress of Chinese Mathematicians. Part 1, 2.} AMS/IP Studies in Advanced Mathematics, vol. 42, pt. 1, 2 (American Mathematical Society, Providence, RI, 2008), pp.\  41--61.

\bibitem[19]{Mori}  \bibaut{S. Mori}, Projective manifolds with ample tangent bundle, {\it Ann. of Math.} {\bf 110} (1979), 593--606.

\bibitem[20]{Muk} \bibaut{S. Mukai}, Biregular classification of Fano 3-folds and Fano manifolds of coindex 3, {\it Proc. Nat. Acad. Sci. USA} {\bf 86} (1989), 3000--3002.

\bibitem[21]{Mu} \bibaut{R. Mu\~ noz}, Varieties with low dimensional dual variety, {\it Manuscripta Math.} {\bf 94} (1997), 427--435.

\bibitem[22]{NO} \bibaut{C. Novelli, G. Occhetta}, Projective manifolds containing a large linear subspace 
with nef normal bundle, {\it Michigan Math. J.} {\bf 60} (2011), 441--462.

\bibitem[23]{OSS} \bibaut{C. Okonek, M. Schneider, H. Spindler}, {\it Vector Bundles on Complex Projective Spaces}, Progress in Mathematics 3, Springer, 1980.

\bibitem[24]{LQEL I}  \bibaut{F. Russo}, Varieties with quadratic entry locus. I, {\it Math. Ann.} {\bf 344} (2009), 597--617.


\bibitem[25]{Ru} \bibaut{F. Russo}, {\it On the Geometry of Some Special Projective Varieties}, Lecture Notes of the Unione Matematica Italiana, Springer, 2016.

%\bibitem[Wa] {Wa} Watanabe

%\bibitem[Sato] {Sato}

\bibitem [26]{Wi} \bibaut{J. Wi\'sniewski},  On a conjecture of Mukai, {\it Manuscripta Math.} {\bf 68} (1990), 135--141.

\bibitem [27]{Zak} \bibaut{F.L. Zak},  {\it Tangents and Secants of Algebraic Varieties}, Transl.  Math. Monogr., vol. 127, Amer. Math. Soc., 1993.

\end{thebibliography}
\end{document}